\documentclass[11pt]{article} 
\usepackage{hfstyle}
\usepackage{graphicx,bbm}
\usepackage[right]{eurosym}
\usepackage{textcomp}
\usepackage[T1]{fontenc}
\usepackage{amsmath}

\usepackage[polutonikogreek,english]{babel}
\usepackage[utf8x]{inputenx}
\newcommand{\greek}[1]{{\selectlanguage{polutonikogreek}#1}}

\usepackage[draft,firstpage]{draftmark}
\draftmarksetup{fontfamily=ptm,angle=0,scale=0.07,mark={\mbox{Published in \emph{Canadian Numismatic Journal}  vol. 61 no. 4 pp. 167-169, 2016 (part 1) and vol. 61 no. 5 pp. 208-213, 2016 (part 2)}},color=red,xcoord=-78,ycoord=120}

\usepackage{hfstyle}
 
\begin{document}
\title{Mathematical formulae on coins}
\author{Henryk Fuk\'s, \#22079 
      \oneaddress{
         Department of Mathematics and Statistics,
 Brock University,\\
     St. Catharines, Ontario, Canada  \\
         \email{hfuks@brocku.ca}
       }
   }

%
\Abstract{
In the article ``The Tale of Two Queens and Two Towering Figures'' published in  CNJ in 2012 (CNJ vol. 57 No. 5, pp. 304-315), we discussed
the contributions of Copernicus and Newton to coin minting and monetary reforms, as well as the commemoration of their
achievements on contemporary coins. In the current article we will examine more explicit aspects of the relationship between
mathematical sciences and numismatics, namely the presence of mathematical formulae on coins. We will survey examples of coins
depicting mathematical formulae, ranging from some universally recognizable ones to more advanced symbolic expressions
which are known only to specialists.
}
\maketitle
As Nicolaus Copernicus once remarked, \textit{mathemata mathematicis scribuntur} - mathematics is written for mathematicians.
This may be so, but some mathematical formulae crossed the boundaries of mathematics and became a part of popular culture. Some of them actually
ended up being reproduced on coins, and these will be a subject of this article.

Although abstract geometric patterns with some apparent mathematical meaning occasionally appeared on ancient coins, the idea to put
a formula on the face of a coin is distinctly modern. In fact, only a handful of late 20th century coins bear inscriptions
utilizing mathematical notation. Such inscriptions are  placed on coins almost exclusively in order to emphasize the achievements
of a mathematician or a physicist who is commemorated by the coin and,  for that reason, the formulae are typically rather secondary elements of the  coin's design. When a mathematician takes a closer look at such coins, some surprises are surely to be found. 
Let us examine some of them.

There is hardly any formula more famous that $E=mc^2$. This expression of equivalence
of mass $m$ and energy $E$, where $c$ denotes the speed of light,  comes from the paper of Albert Einstein ``Ist die Trägheit eines Körpers von seinem Energieinhalt abhängig?'' (``Does the inertia of a body depend upon its energy-content?'')  published in \textit{Annalen der Physik} in 1905 \cite{Einstein1905}. In the popular mind,
$E=mc^2$ functions as an integral attribute of Einstein, often serving as a symbol of the theory of relativity, or, even more
often, as an exemplification of Einstein's genius. It is not surprising, therefore, that it is prominently featured on several
commemorative coins. The oldest one that I am  aware of  is  the Chinese 10 yuan silver issue from 1991. It shows the formula 
in the background of the portrait of Einstein, who appears as if he had just written it on the blackboard behind him.
German 5 euro coin issued in 2005, commemorating the centennial anniversary of the special theory of relativity, features
$E=mc^2$ in the center of the reverse. In the same year, designated as ``Einstein Year'' or ``The Year of Physics'', Israel  issued a  1 shekel coin displaying the formula along the bottom border of the reverse, with Einstein's autograph reproduced on the obverse. An interesting feature of this coin is that the circular lines on the reverse, if viewed from a certain angle, create the image of Albert Einstein. 
\begin{center}
  \includegraphics[width=4.5cm]{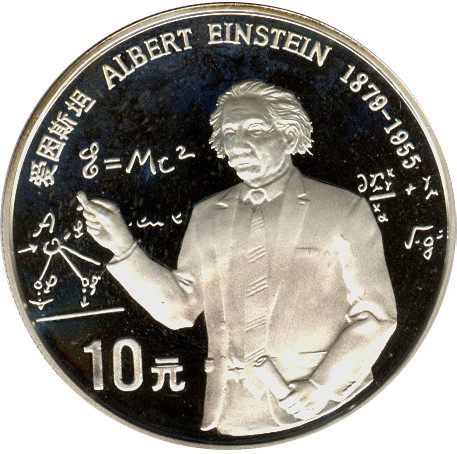}  \includegraphics[width=4.5cm]{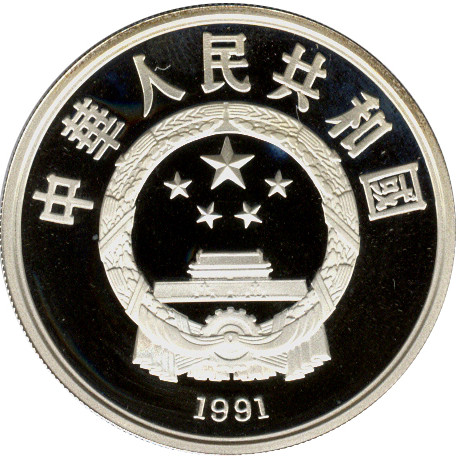} 
 \end{center}
\begin{center}
  \includegraphics[width=4.5cm]{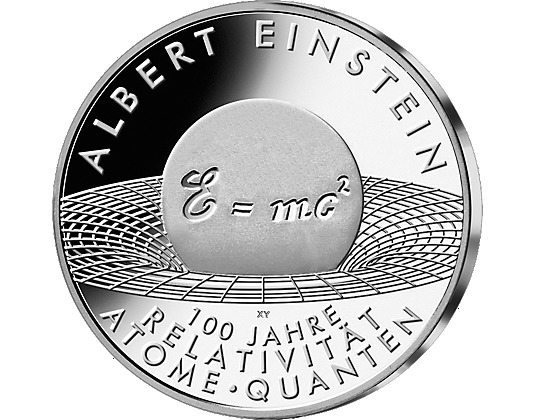} \includegraphics[width=4.5cm]{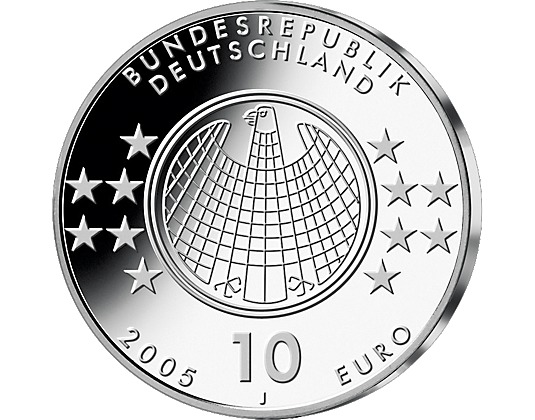}
 \end{center}
\begin{center}
  \includegraphics[width=4.5cm]{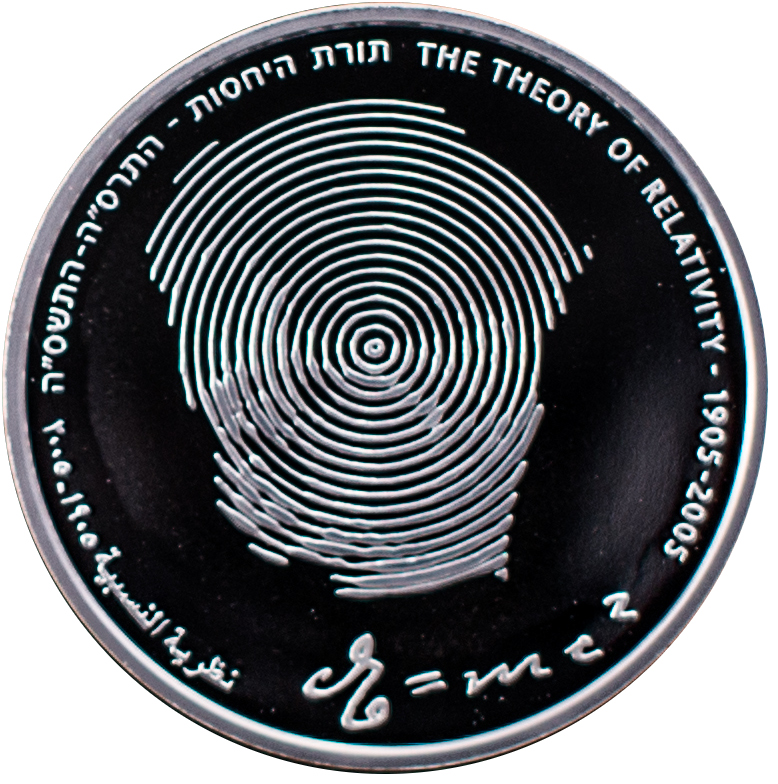}  \includegraphics[width=4.5cm]{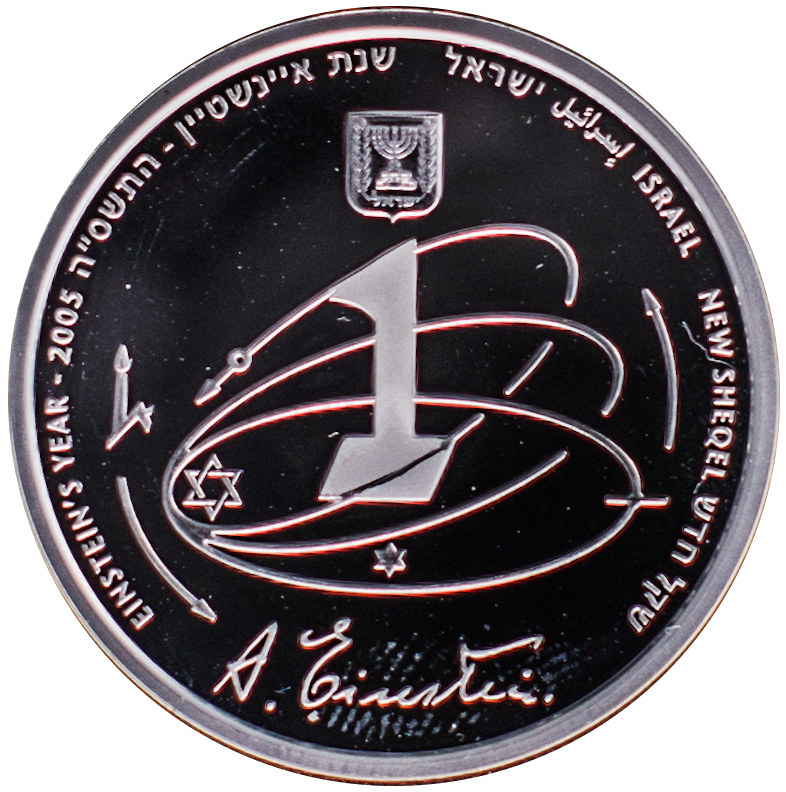} 
 \end{center}
In 2015, on the 100th anniversary of general theory of relativity, the Royal Canadian Mint produced a giant 10 oz. silver coin 
with \$100 denomination, also featuring $E=mc^2$.  The portrait of Einstein on this coin is the work of   
   Canadian photographer Yousuf Karsh (1908-2002). In addition,  silver and gold  bullion coins with $E=mc^2$ privy mark below the maple leaf were produced in the same year. 
 \begin{center}
  \includegraphics[width=4.5cm]{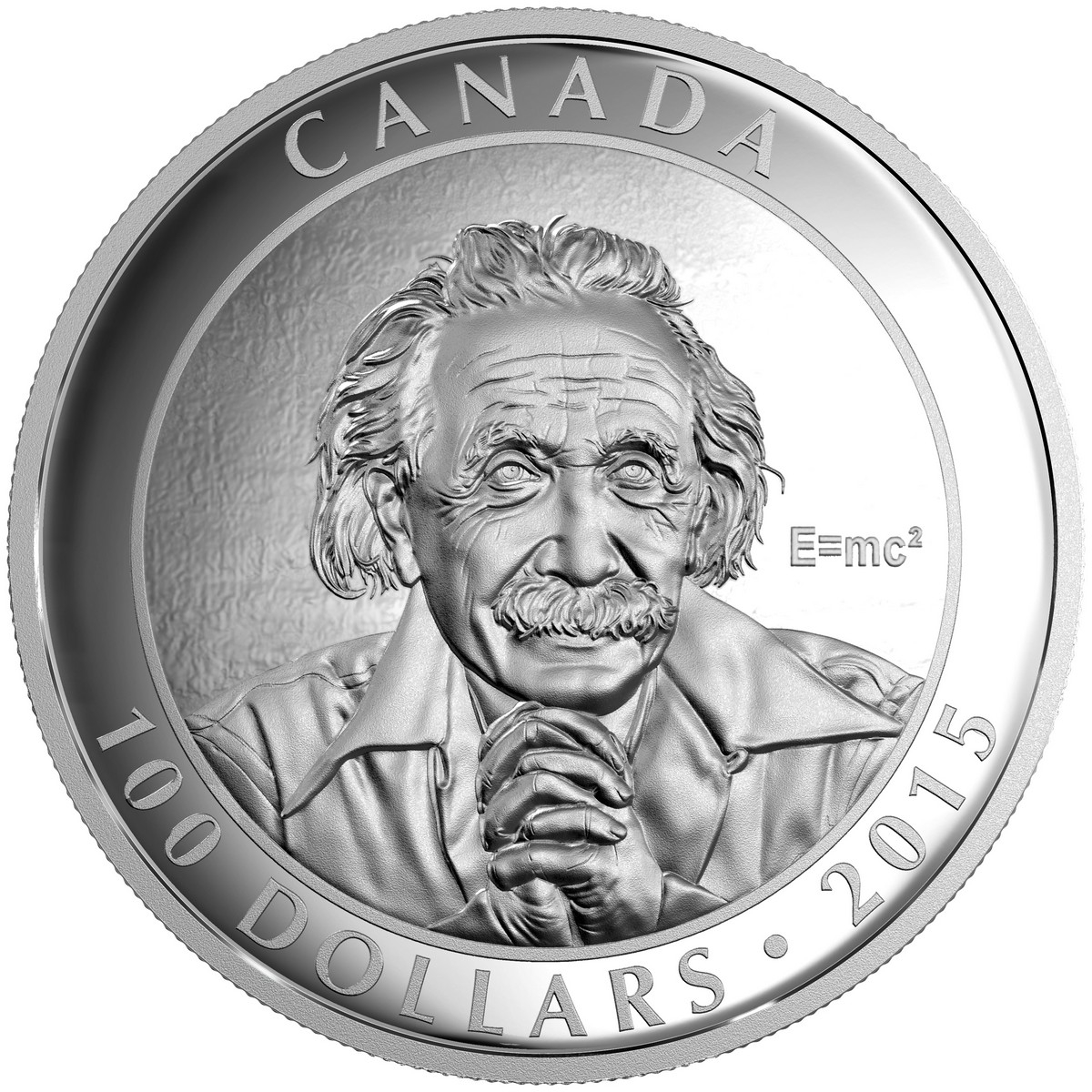}  \includegraphics[width=4.5cm]{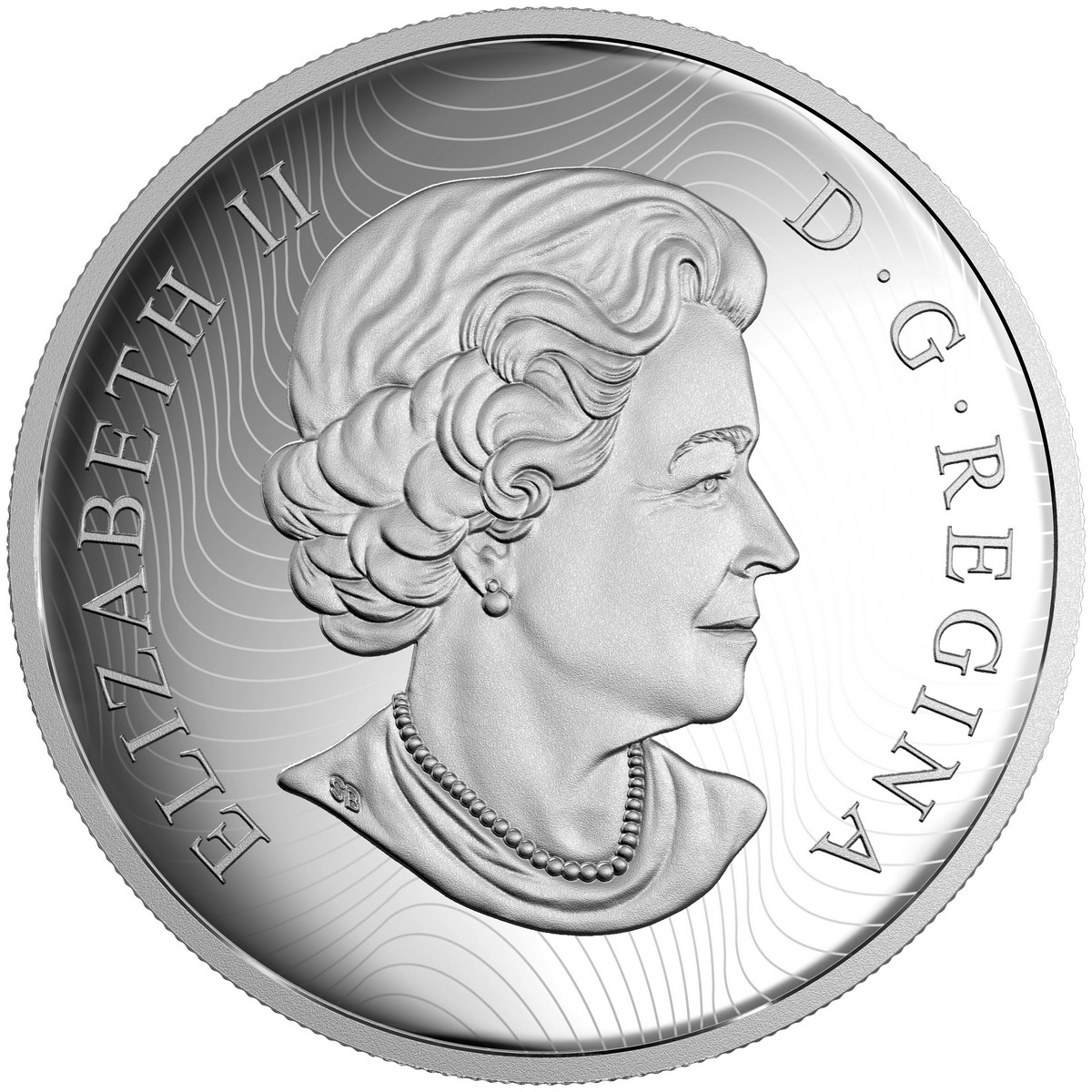} 
\includegraphics[width=4.5cm]{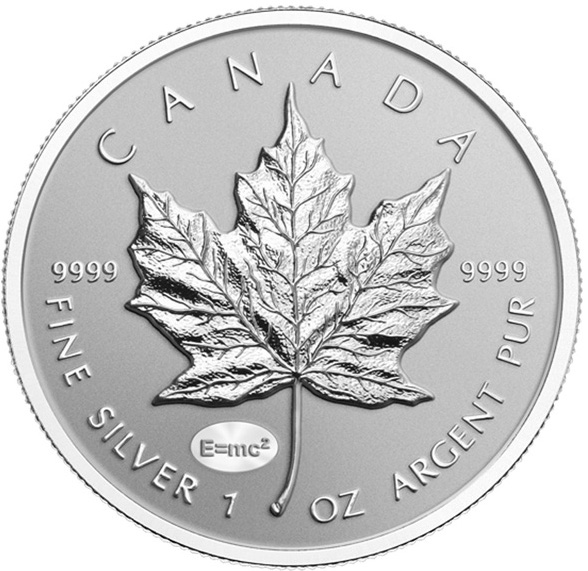}
 \end{center} 

The energy equivalent to mass $m$ given by the equation  $E=mc^2$ is sometimes called the ``rest energy'', because
it is an internal or intrinsic energy of a stationary body. If an object is moving with velocity $v$, it also possesses some kinetic energy, and in this  case its total energy is  given by
\begin{equation*}
 E=\frac{m c^2} {\sqrt{1-v^2/c^2}}.
\end{equation*}
A stylized version of this more general energy formula appeared on the reverse of the bimetallic 500 Lire coin of San Marino, which was part of
the 1984 thematic series ``Science for People''. The core of the coin, made of bronzital, features a portrait
of Einstein, while the  formula is engraved on the  silver-colored outside ring (made of acmonital).
   Canadian photographer Yousuf Karsh (1908-2002). In addition,  silver and gold  bullion coins with $E=mc^2$ privy mark below the maple leaf were produced in the same year. 
 \begin{center}
  \includegraphics[width=4.5cm]{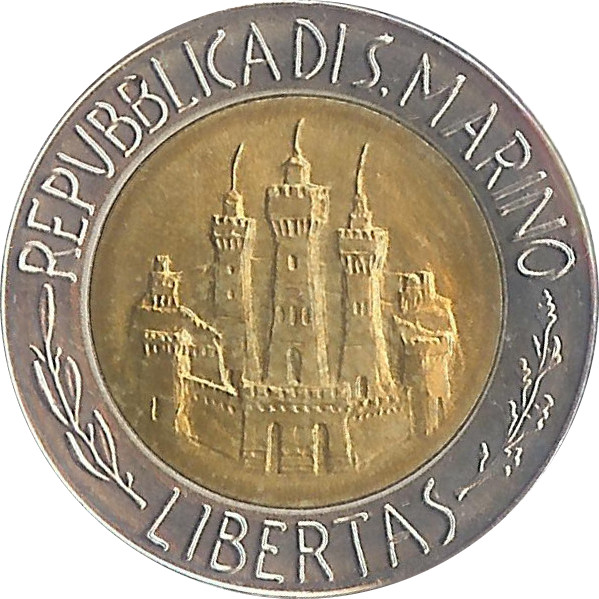}  \includegraphics[width=4.5cm]{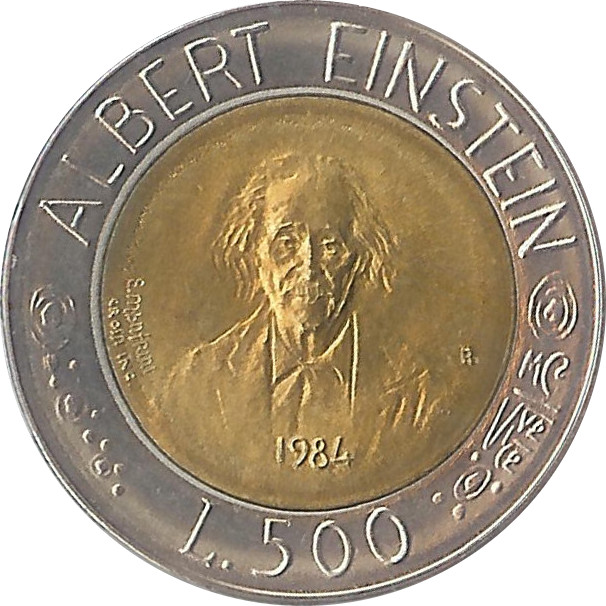} 
\includegraphics[width=2.5cm]{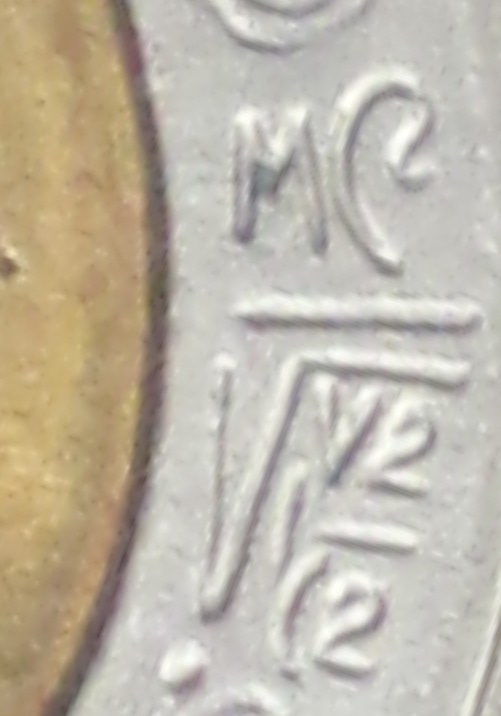}
 \end{center} 

The equivalence of mass and energy equation is not the only equation connected to relativity and Einstein appearing on coins.
In 1979 Switzerland minted a pair of coins to commemorate the 100th anniversary of Einstein's birth, one with his portrait 
and another one with two formulae engraved on the reverse. The first of these formulae is $R_{uv}=0$, which indicates components of Ricci curvature tensor being zero. In differential geometry, the Ricci tensor is used to ``measure'' how much the space geometry differs from Euclidean
geometry. In the general theory of relativity, Ricci tensor is zero in empty space.  It is not clear to me why did the designer of the coin
select this particular formula, since it describes a rather special case of vacuum. It seems that central equations of general relativity,
known as Einstein filed equations, would be more appropriate. These describe the relationship between the Ricci tensor and the stress–energy tensor 
and are normally written as follows,
\begin{equation*}
 R_{\mu \nu} - \frac{1}{2}R \, g_{\mu \nu} + \Lambda g_{\mu \nu}= \frac{8 \pi G }{c^4} T_{\mu \nu}.
\end{equation*}
One can only guess that the full form of field equations was considered too long, too complex, or not fitting to the artistic concept of the coin.
The second equation engraved on the coin, $\delta \int ds=0$, represents the so-called principle of least action, written in abbreviated form.
Without going into details, it is a principle used to derive equations of motion of mechanical systems, based
on a branch of mathematical analysis know as calculus of variations. Einstein field equations can be derived  through the principle of least action when applied to Einstein–Hilbert action, and this is probably why the principle has been selected to appear on the coin.
\begin{center}
  \includegraphics[width=4.5cm]{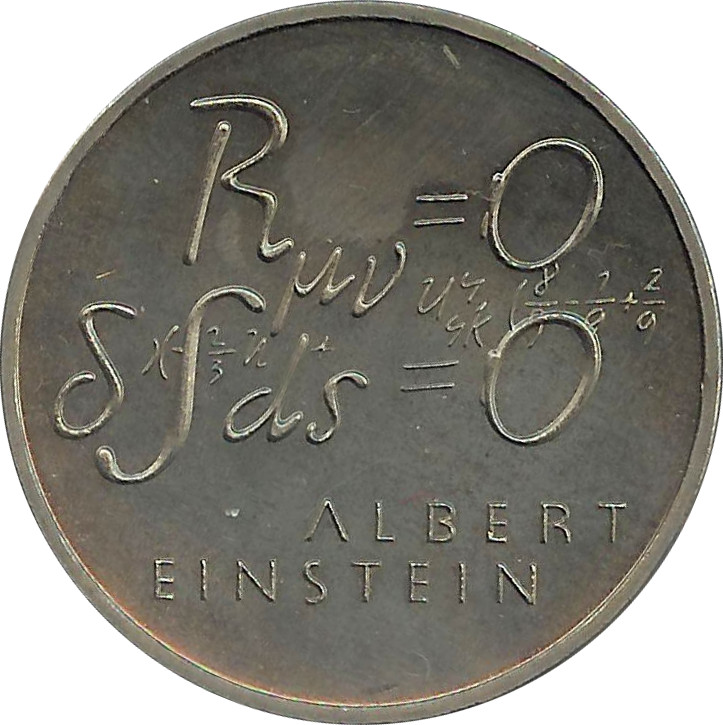}  \includegraphics[width=4.5cm]{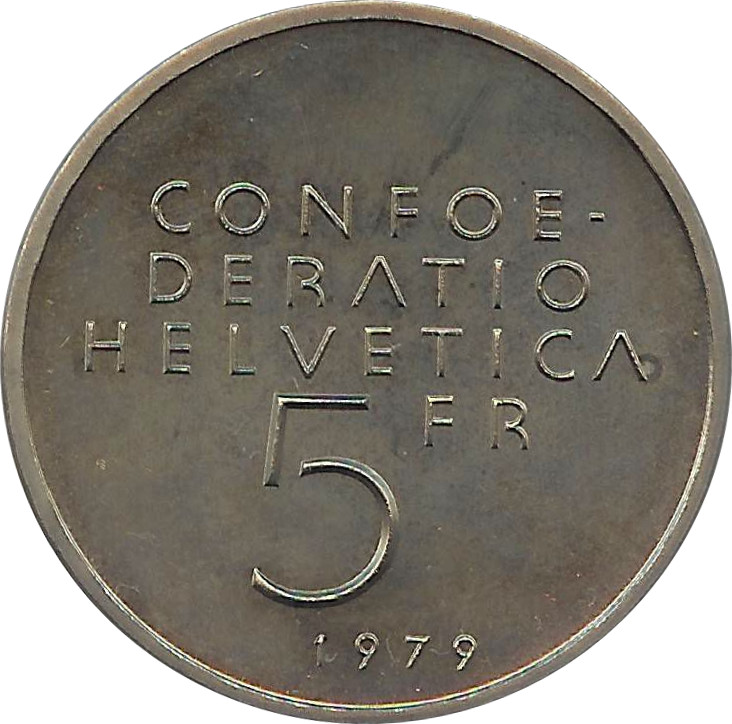} 
 \end{center}

Speaking of Einstein, one cannot fail to mention that among all scientists, his portrait appears on coins most often. In my collection I have 
more than 10 coins featuring Einstein, issued by countries as diverse as Germany, Switzerland, San Marino, Israel, China, Paraguay, Uganda,
 Togo, and others. 
New coins with Einstein appear almost every year, thus one cannot help but ask: is Einstein's contribution to science really
commensurate with his fame and popularity? Not every physicist would consent. As important as relativity theory is in physics,
its impact on our our lives is rather minimal. Other achievements of 20th century physics, in particular the development of quantum mechanics,
were far more important and led to the rapid advances in technology, yet one does not see too many coins commemorating
quantum mechanics. Moreover, Einstein's level of contribution to development of special theory relativity was always a subject of bitter disputes 
in the physics community. It is certain that Henri Poincar\'e (1954--1912) laid foundations to the special theory of relativity well before Einstein, and that Einstein
was familiar with his work. Some recent historians go as far as accusing Einstein of plagiarism of Poincaré's results. Jules
Levegue in his 2004 book \cite{Lev2004} presents quite convincing evidence of this
claim. In the meanwhile, nobody refers to $E=m c^2$ as ``Poincar\'e-Einstein formula'', and Henri Poincar\'e, one of the greatest mathematicians of all times,
has yet to see to a coin commemorating his achievements.

Having discussed the most popular of all formulae, let us turn our attention to other entrants. We will present them in
chronological order, starting from the oldest one, the Pythagorean equation. The triangular ``millennium coin of Uganda displays
the formula $\alpha^2+\beta^2=\gamma^2$ below the portrait of Pythagoras (c. 570 -- c. 495 BC), credited with the first proof
of the theorem commonly known as the Pythagorean theorem. 
\begin{center}
  \includegraphics[width=6.5cm]{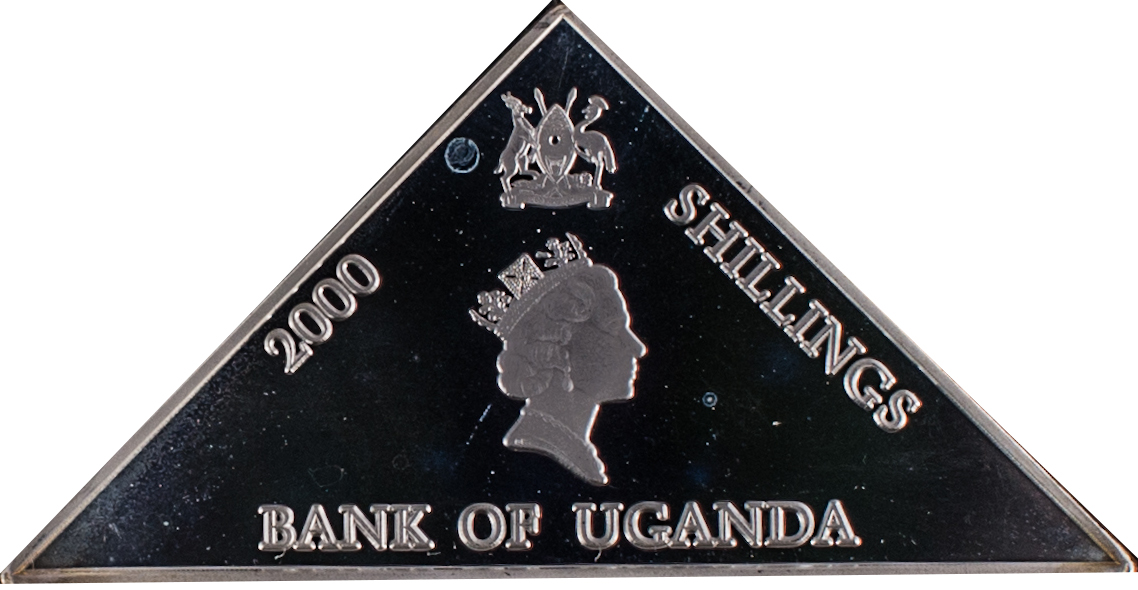}  \includegraphics[width=6.5cm]{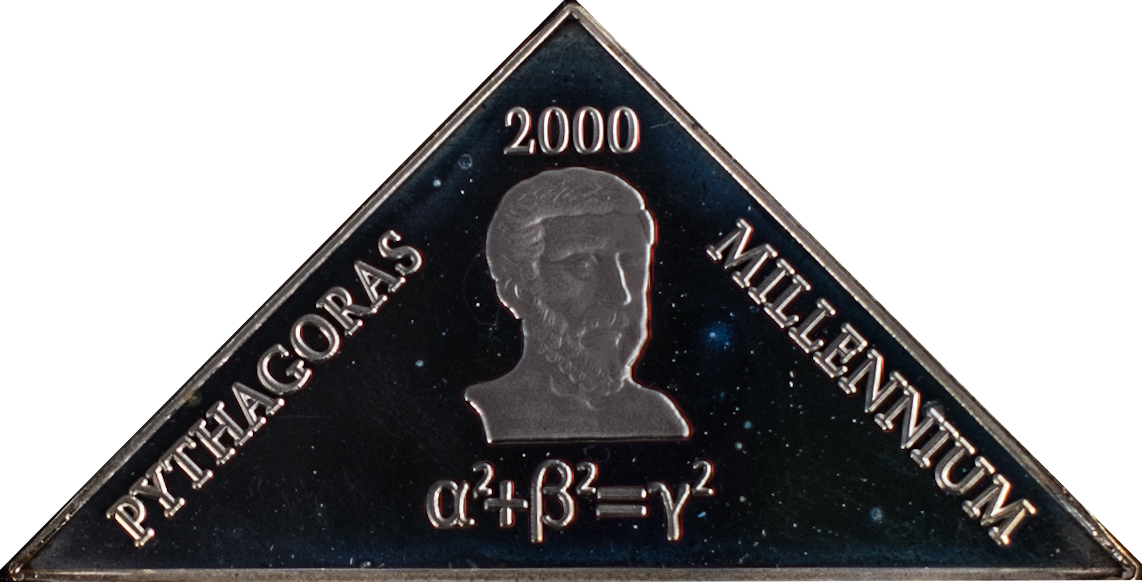} 
 \end{center}
The theorem states that in a right triangle, the square of the hypotenuse  is equal to the sum of  squares of the other two sides. The equation on the coin, written in Greek letters, is more commonly
represented using Latin letters as $a^2+b^2=c^2$, where $a$, $b$, and $c$ are the lengths of the sides of the right triangle.
It should be stressed, however, that even if written in Greek letters, the formula $\alpha^2+\beta^2=\gamma^2$ uses modern
algebraic notation, unknown to ancient Greeks. In book I of Euclid's ``Elements'', written around 300 BC, the Pythagorean theorem is 
stated as proposition 47, as follows:
\begin{quotation}
\greek{Εν τοῖς ὀρθογωνίοις τριγώνοις τὸ ἀπὸ τῆς τὴν ὀρθὴν γωνίαν ὑποτεινούσης πλευρᾶς τετράγωνον ἴσον ἐστὶ τοῖς ἀπὸ τῶν τὴν ὀρθὴν γωνίαν περιεχουσῶν πλευρῶν τετραγώνοις.}
\end{quotation}
In English translation this means ``In right-angled triangles the square on the side subtending the right angle is equal to the squares on the sides containing the right angle''. Of course, it would not be practical to put such a long inscription on the coin, thus we 
will not hold it against the designer that he opted for a modern formula instead.

The next formula  is perhaps less known by the general public, but every mathematician recognizes it instantly. The infinite sum
\begin{equation*}
 \sum_{i=1}^\infty \frac{1}{n^2}=\frac{\pi^2}{6}
\end{equation*}
appears on the Russian 2 ruble coin issued in 2007. 
The origins of this sum go back to the book of the Italian mathematician and  Catholic priest Pietro Mengoli (1626-1686) 
``Novae quadraturae arithmeticae, seu de additione fractionum'',  published in Bologna in 1650.
\begin{center}
  \includegraphics[width=4.5cm]{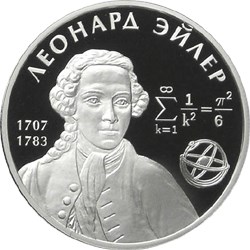}  \includegraphics[width=4.5cm]{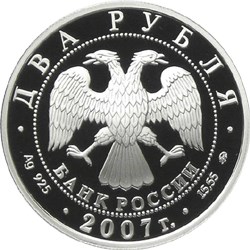}\\  \includegraphics[height=6cm]{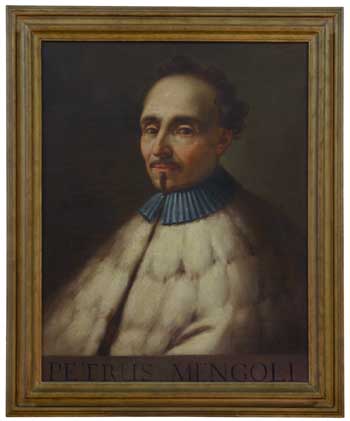}  \includegraphics[height=6cm]{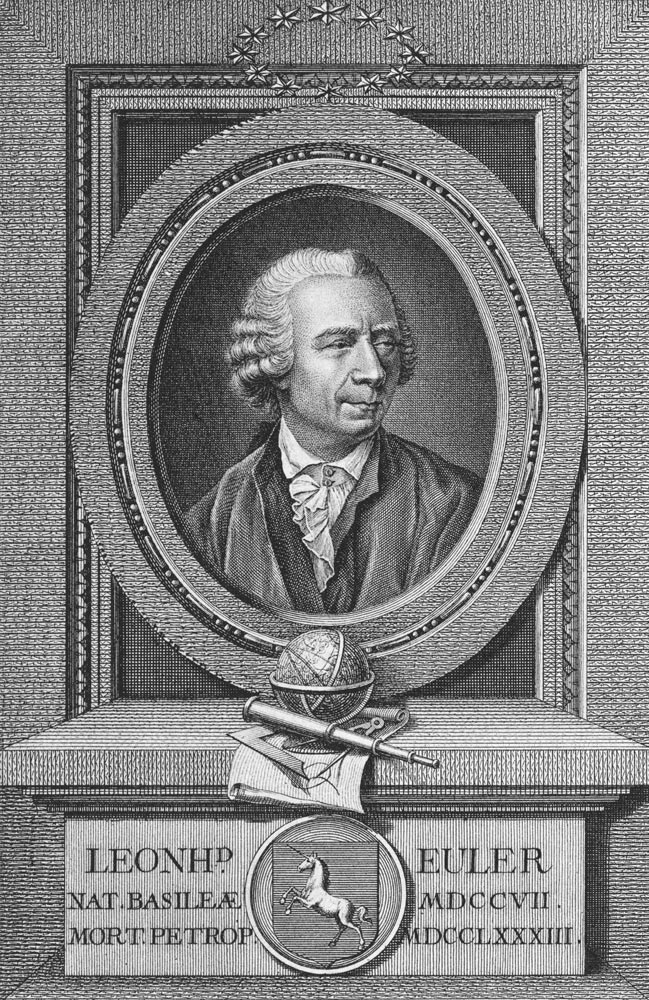} 
 \end{center}
In this work, among other things, Mengoli computed  sums of the reciprocals of various subsets of integers.
There was, however, one sum which he could not compute and left it as an open problem, namely the summation of the reciprocals of the squares of the natural numbers.
The problem remained unsolved until 1734, when Leonard Euler (1707-1783), then twenty eight years old, found that
\begin{equation*}
 \frac{1}{1^2}+\frac{1}{2^2}+\frac{1}{3^2}+\ldots =\frac{\pi^2}{6}.
\end{equation*}
The solution to this problem, sometimes known as the Basel problem after the hometown of Euler, brought him  instant fame
and recognition by the the mathematics community of the day. One should add that the use of the Greek letter sigma to denote summation,
as it appears on the coin, is also due to Euler, who introduced it in 1755 in his book ``Institutiones calculi differentialis''.
In the original publication presenting the solution of the Basil problem \cite{Euler1740}, the summation symbol is not used,
instead, Euler wrote
\begin{equation*}
 \mathrm{I} + \frac{\mathrm{I}}{2^2} +  \frac{\mathrm{I}}{3^2} +  \frac{\mathrm{I}}{4^2} + \frac{\mathrm{I}}{5^2} \mathrm{\,\,\,etc.}=\frac{p^2}{6},
\end{equation*}
using, as it was the custom, I to denote 1 and $p$ to denote $\pi$.

Euler, the most prolific mathematician of all times and certainly one of the most influential ones, resided in St. Petersburg in years
1726--1741 and again from 1766 till his death. He was also a member of the Saint Petersburg Academy of Sciences, thus it is not surprising
that the Moscow Mint decided to honor the 300-th anniversary of his birth with a silver proof coin. The decision to incorporate the
Basel problem into the design of the coin was, in my opinion, a very good choice. The solution of this problem was indeed
a turning point in Euler's career, and his proof that the sum is equal to $\pi^2/6$ is considered a typical ``Eulerian'' proof,
a mark of true genius. Nevertheless, given that Euler wrote over 800 research papers in his lifetime, I am sure that it was not easy
to choose just one formula to represent his vast legacy.

One of the leading mathematicians of 19-th century Russia, often considered to be an intellectual descendant of Euler, was
Mikhail Vasilyevich Ostrogradsky (1801--1862). In 2001, to commemorate the 200-th anniversary of his birth,
the mint of the National Bank of Ukraine issued 2 a hryvnia coin, the first in the series ``Outstanding Personalities of Ukraine''.
The coin, in addition to Ostrogradsky's portrait, features the integral $\int\frac{P(x)}{Q(x)}dx$. This is, unfortunately,
not a complete formula, but rather a part of a formula known as the Ostrogradsky method of integrating rational functions \cite{Ostro1845}.
\begin{center}
  \includegraphics[width=4.5cm]{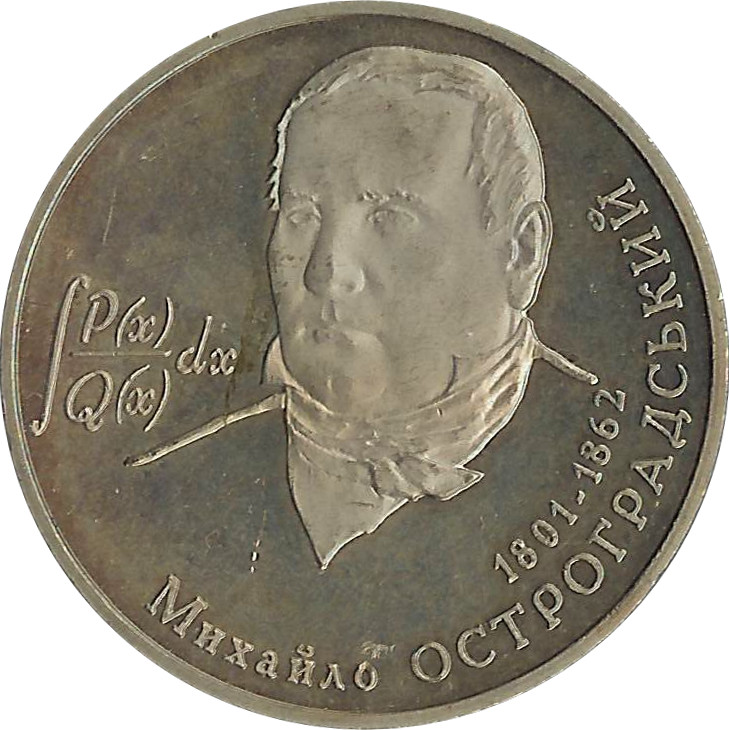}  \includegraphics[width=4.5cm]{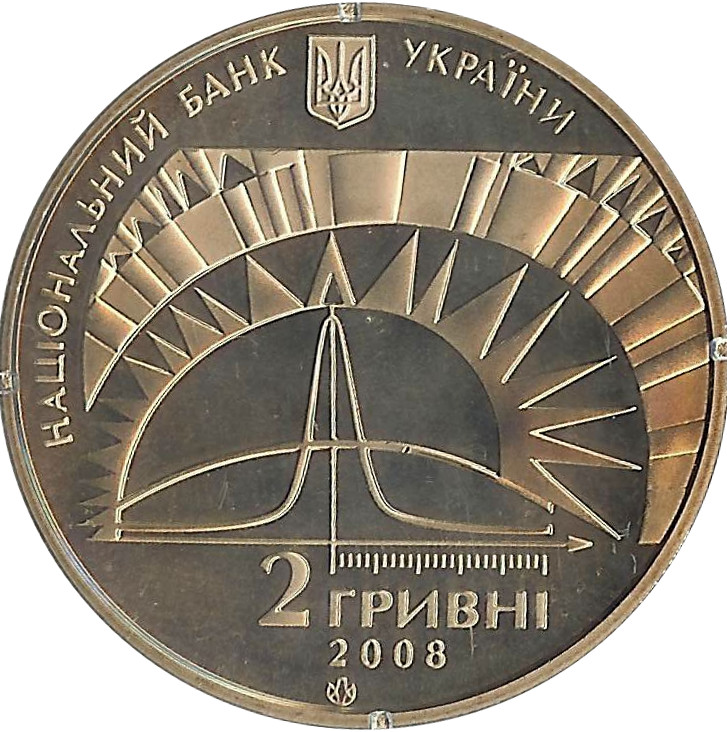}
 \end{center}
Its complete statement is as follows. If $P(x)$ and $Q(x)$ are polynomials with real coefficients, $P(x)/Q(x)$ is a proper fraction, 
and $Q(x)$ has multiple roots, then
\begin{equation*}
 \int\frac{P(x)}{Q(x)}dx= \frac{P_1(x)}{Q_1(x)} + \int \frac{P_2(x)}{Q_2(x)},
\end{equation*}
where $Q_1(x)$ is the greatest common divisor of $Q(x)$ and $Q^{\prime}(x)$, and $Q_2(x)=Q(x)/Q_1(x)$. The polynomials $P_1$ and
$P_2$ are obtained by matching corresponding coefficients on both sides of the equation, after differentiating it. What we
see on the coin is just the left hand side of the above equation. As in the case of Einstein's field equations, designers
presumably wanted to avoid cluttering the coin, thus deciding to include only a truncated part.

The next scientist honored by a coin with a formula is not very well know today, especially in the West. 
Pyotr Nikolaevich Lebedev (1866--1812) was a Russian physicist who was the first to measure the pressure of light in 1899 \cite{Leb1901}.
The 1 ruble coin celebrating the 125-th anniversary of his birth is one of the last commemorative coins issued by the Soviet Union before its formal dissolution in 
December 1991. The inscription on the coin is the expression of the light pressure $p$ in terms of the density of photons $N$,
frequency $\nu$, and the coefficient of reflection $\rho$,
\begin{equation*}
 p=N \frac{h \nu}{c}(1+\rho),
\end{equation*}
where $h$ is the Planck constant and $c$ is the speed of light. This was actually a circulating coin, with a total mintage of over 2 million.
\begin{center}
  \includegraphics[width=4.5cm]{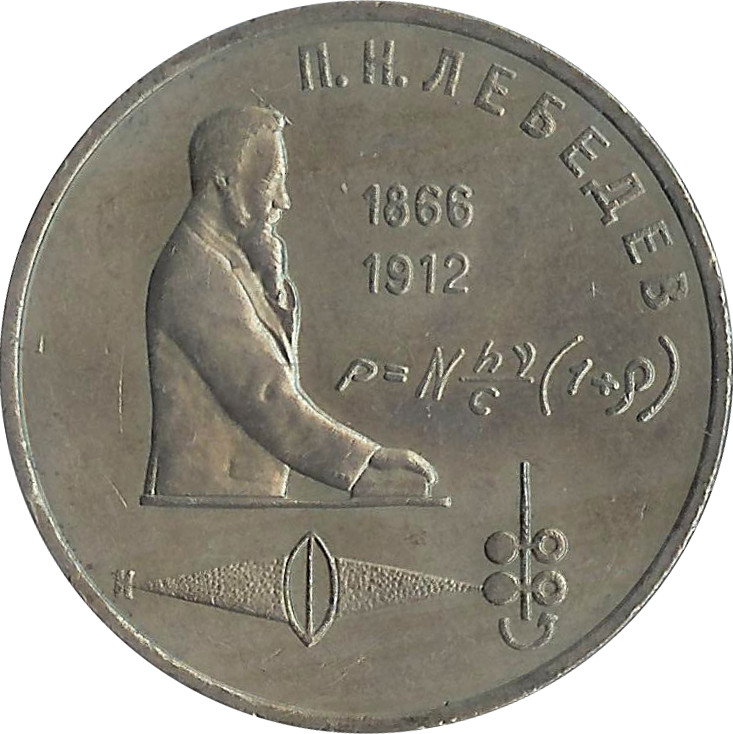}  \includegraphics[width=4.5cm]{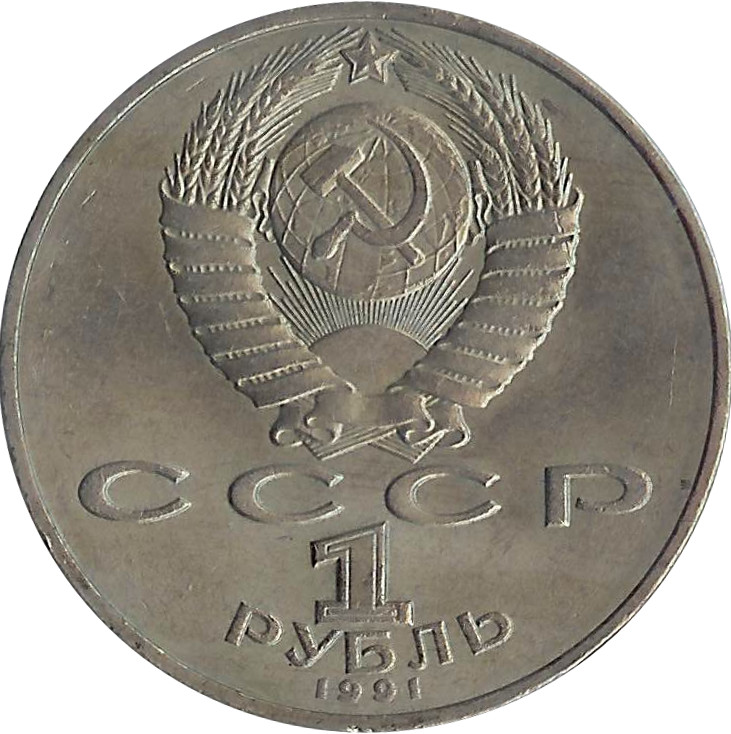}
 \end{center}

In the same series of coins as the issue honoring Ostrogradsky described earlier, in 2009 Ukraine commemorated the centennial anniversary of the birth of a famous mathematician and theoretical
physicist of the Soviet era, Nikolay Nikolayevich Bogolyubov (1909--1992), one of the founders on nonlinear
mechanics and a significant contributor to the quantum field theory and statistical physics. This 2 hrywni Ni-Ag alloy
coin features interesting formulae on both sides. Both of them belong to rather advanced 20th century mathematics,
thus we will not be able to explain their full meaning here. 
\begin{center}
  \includegraphics[width=4.5cm]{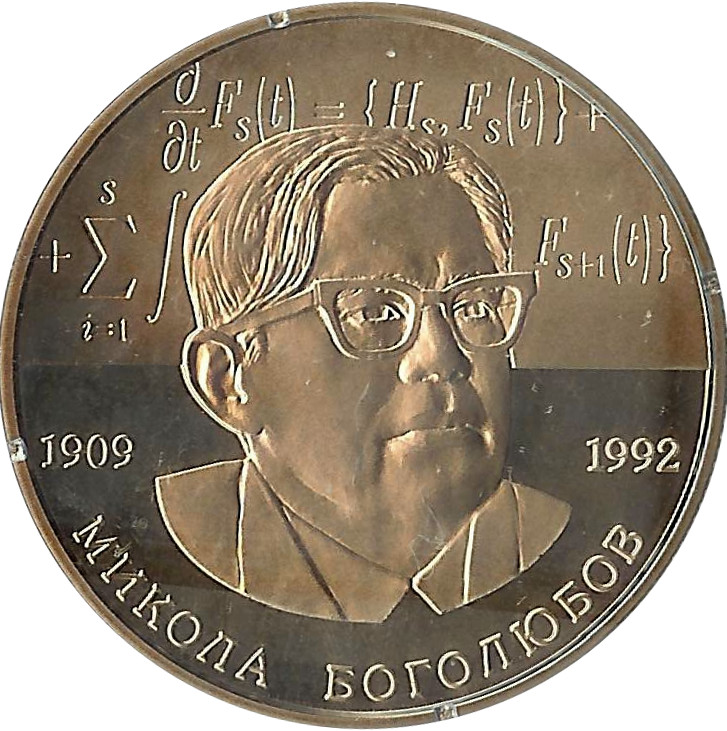}  \includegraphics[width=4.5cm]{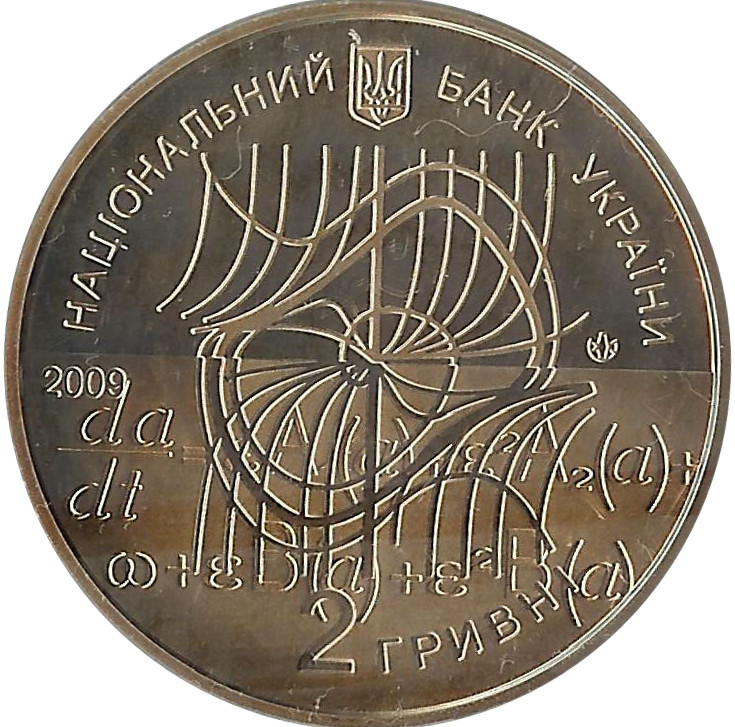}
 \end{center}

The formula on the reverse is today known as the Bogoliubov–Born–Green–Kirkwood–Yvon hierarchy, 
sometimes called simply the Bogoliubov hierarchy \cite{Bogo46}. It is  a set of equations describing the dynamics of a system of a large number of interacting particles. The equation for an $s$-th particle distribution function $F_s$ depends on the $(s + 1)$-th particle distribution function $F_{s+1}$, thus forming  a  hierarchy of coupled equations, 
\begin{equation*}
 \frac{\partial}{\partial t} F_s(t)=\{H_s,F_s(t)\}+\sum_{i=1}^{s} \int \frac{1}{v}\{ \phi(q_i-q), F_{s+1}(t)\} dx.
\end{equation*}
Note that on the coin, part of the expression under the integral sign is obscured by the Bogoliubov's portrait.

The obverse features two formulae, partially covered by a graph which we will describe shortly. These formulae
are
\begin{equation*}
 \frac{da}{dt}=\epsilon A_1(a)+\epsilon^2 A_2(a)+ \ldots, \mbox{\,\,\,\,\,and\,\,\,\,\,}
  \frac{d \psi}{dt}=\omega + \epsilon B_1(a)+\epsilon^2 B_2(a)+ \ldots.
\end{equation*}
They belong to the Krylow-Bogoliubov-Mitropolsky method \cite{mickens1981} of solving differential equations
describing nonlinear osillations of the type
\begin{equation*} 
 \frac{d^2y}{dt^2}+ \omega y =\epsilon F(y, \dot{y}),
\end{equation*}
where one seeks the solution in the form 
\begin{equation*}
 y=a \cos \phi + \epsilon u_1(a,\psi)+\epsilon^2 u_2(a,\psi)+ \ldots.
\end{equation*}
The quantities $a$ and $\psi$ must satisfy the two differential equations shown on the coin. The graph partially covering the
equations is a phase portrait (graphical representation of the solution) of some nonlinear oscillator, most likely van der Pol oscillator, which is a well known
example of a nonlinear oscillator.

Speaking of the 20th century mathematics, one of its most prominent representatives has been honored not with one, but with three
coins.
In 2012, the National Bank of Poland released into circulation three coins commemorating the great Polish mathematician
Stefan Banach (1892--1945). Banach was  the principal  founder of the branch of mathematics known as the functional analysis
and  also one of the founders of Polish School of Mathematics. The coins are made of gold (200 zł), silver (10 zł) and 
Nordic Gold alloy (2 zł). Each of them depicts a formula named after Banach.
\begin{center}
  \includegraphics[width=4.5cm]{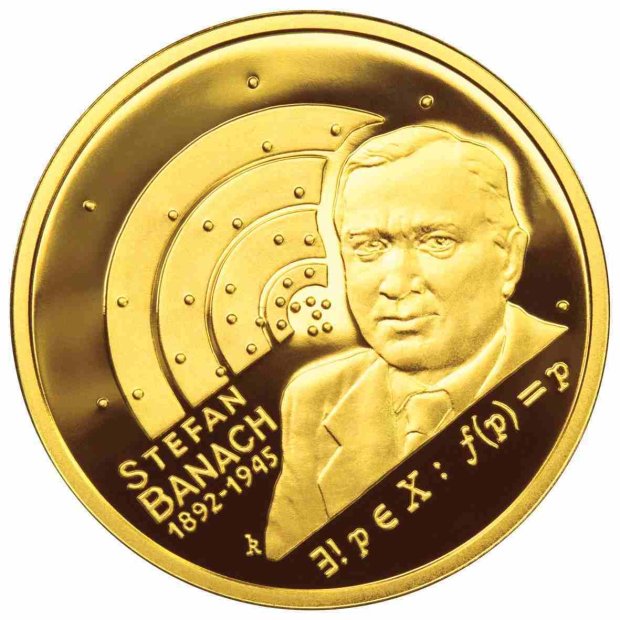}  \includegraphics[width=4.5cm]{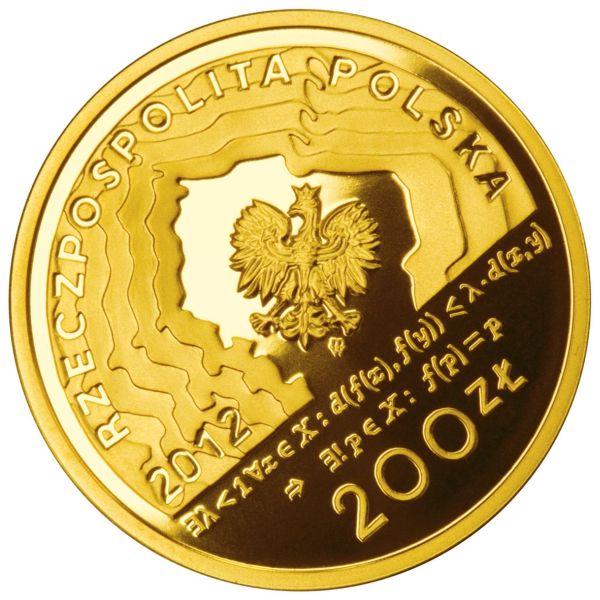} 
 \end{center}
The obverse of the golden coin features  the Banach fixed-point theorem, also known as the contraction mapping theorem,
proved by Banach and published in his seminal paper  ``Sur les 
op\'erations dans les ensembles 
abstraits et leur 
application aux \'equations
int\'egrales'' \cite{Banach1922}.

It says that  every non-empty complete metric space with a contraction mapping has a unique fixed point. This statement is
written on the coin in the symbolic form,
\begin{equation*}
 \exists \lambda < 1 \,\, \forall x,y \in X: d(f(x), f(y)) \leq \lambda \cdot d(x,y)
 \Rightarrow  \exists ! p \in X: f(p)=p.
\end{equation*}
Unfortunately, a close inspection of the coin reveals that the engraving of the formula contains an error: $\forall x,y \in X$
is written as $\forall x \in X$, with $y$ omitted. 
\begin{center}
\includegraphics[width=4.5cm]{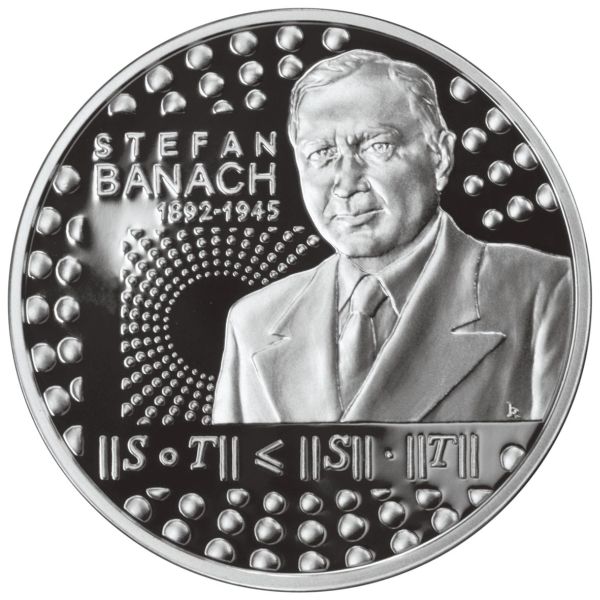}  \includegraphics[width=4.5cm]{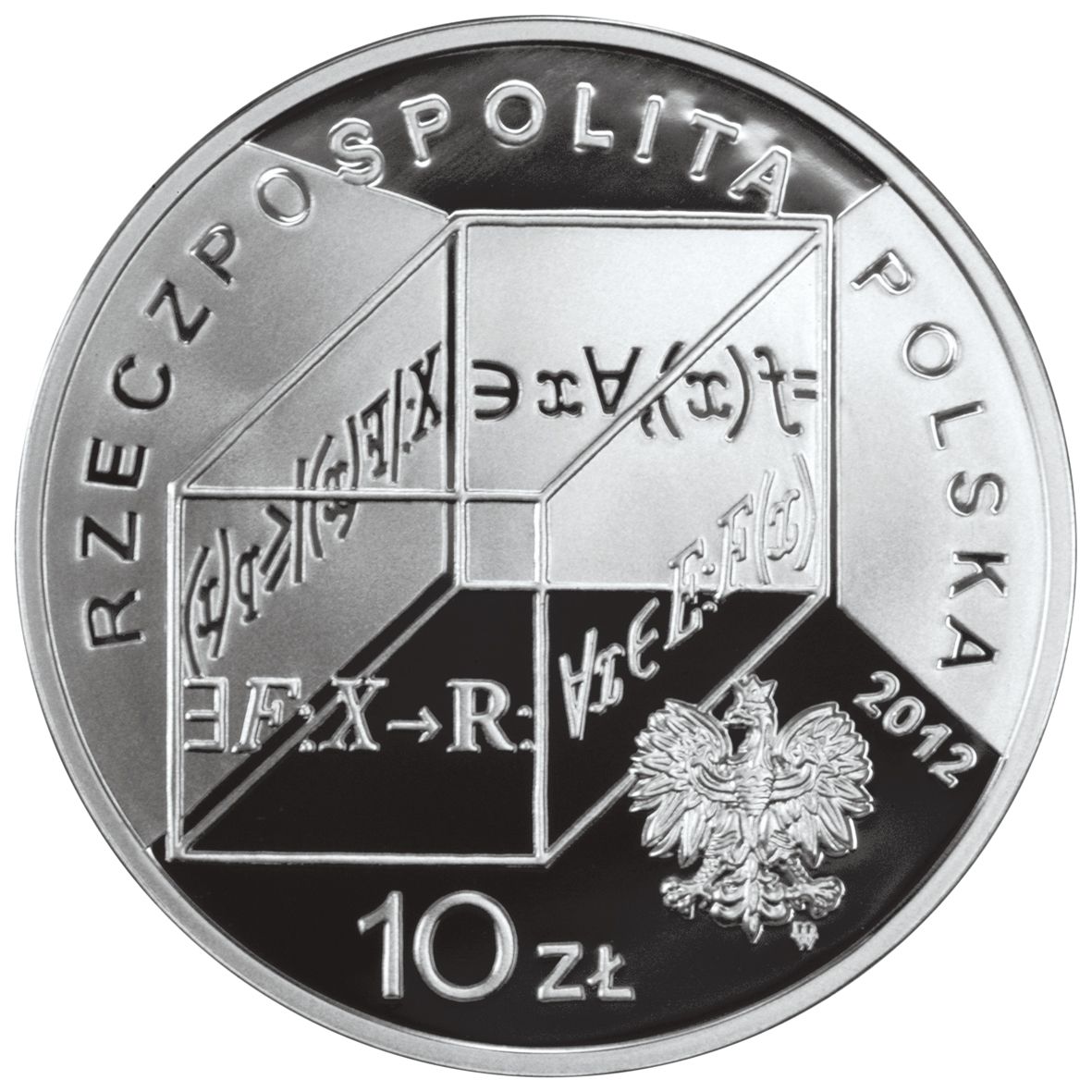}
 \end{center}
The reverse of the silver coin shows, under the portrait of Banach, the formula $||S \circ T|| \leq ||S|| \cdot ||T||$.
In functional analysis  one proves that this inequality is  satisfied by every pair of bounded operators in a Banach space.
It is also used in the definition of Banach algebras.

The obverse of the silver coin features the following sequence of symbols inscribed on walls of a cube,
\begin{equation}
 \exists F:X \to \mathbbm{R}: \forall x \in E: F(x) = f(x), \forall x \in X: |F(x)| \leq p(x).
\end{equation}
This is a conclusion (hypothesis is omitted) of a central theorem of functional analysis, known and the Hahn-Banach theorem, independently
proved by Banach \cite{Banach1922} and an Austrian mathematician Hans Hahn (1879--1934). Its meaning is rather complicated to explain in a short space,
thus we will be satisfied by saying that it allows to extend bounded linear functionals from a subspace to the entire space.
\begin{center}
\includegraphics[width=4.5cm]{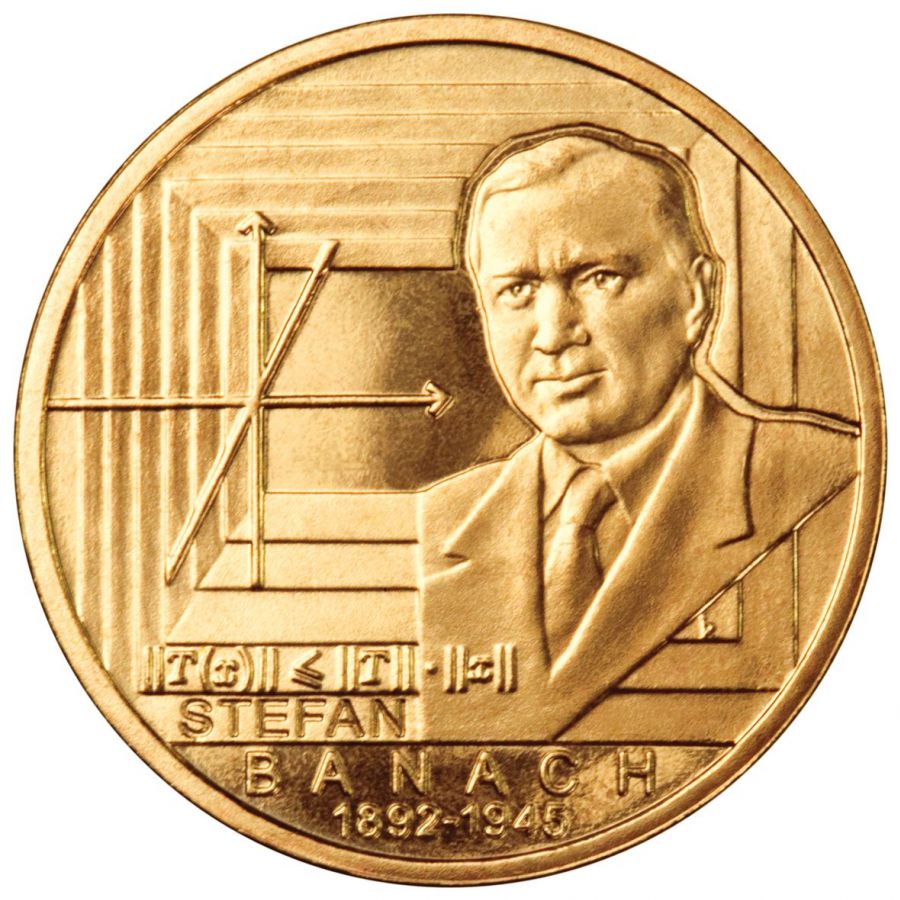}  \includegraphics[width=4.5cm]{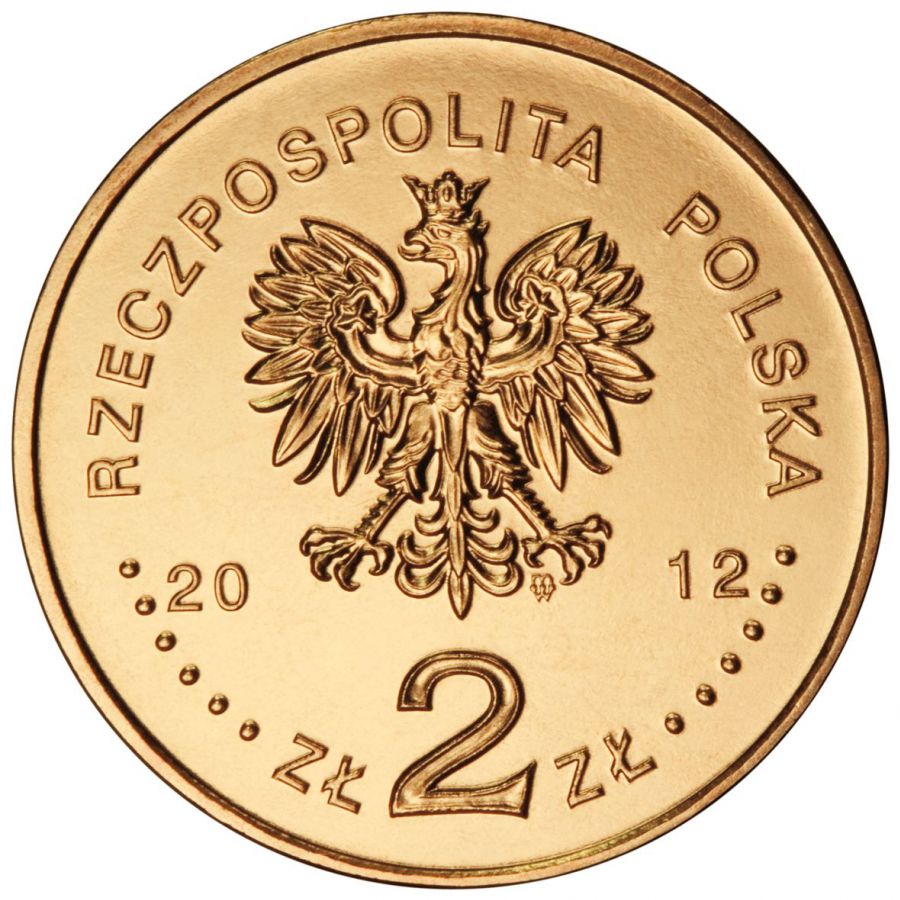}
 \end{center}
Finally, the reverse of the 2 zł coin depicts, under the portrait of Banach, the inequality $||T(x)|| \leq ||T|| \cdot ||x||$.
This is a condition satisfied by linear operators in Banach spaces, and is a direct consequence of the definition of the operator norm. 

We will finish on a slightly lighter note. In 2005, which as we recall was celebrated as the World Year of Physics, 
the Central Bank of Ireland produced a beautifully crafted  silver proof 10 euro  coin
to mark the 200th anniversary of the birth of William Rowan Hamilton (1805 – 1865), renowned as one of the world's greatest mathematical
scientists. 
\begin{center}
  \includegraphics[width=4.5cm]{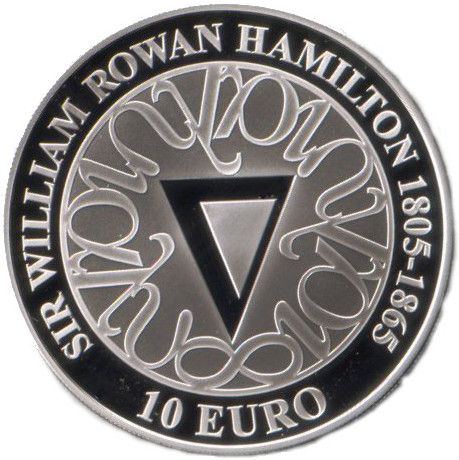}  \includegraphics[width=4.5cm]{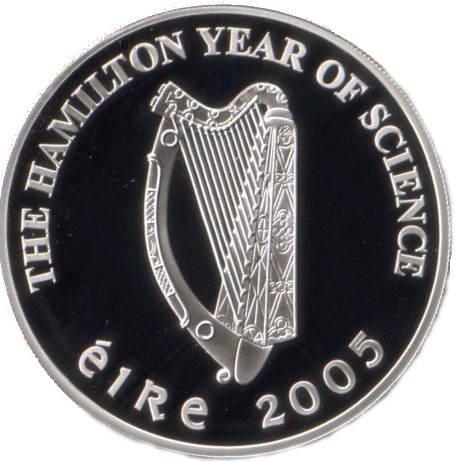}
 \end{center}
The center of the obverse is occupied not by a formula, but by a mathematical symbol which is a part
of great many formulae. The symbol, in the form of inverted Greek capital letter delta, is called the del or nabla operator,
named after the Assyrian harp whose shape it resembles. Of course, the  harp also appears in the coat of arms of Ireland shown
on the reverse,  thus both sides of the coin nicely harmonize.  The nabla operator, defined as 
$\nabla = (\frac{\partial}{\partial x},\frac{\partial}{\partial y},\frac{\partial}{\partial z})$, has been invented by Hamilton,
but the first documented use of this symbol comes from the book of Peter G. Tait (1831-1901) ``An Elementary Treatise on Quaternions'', 
written 1867. Tait calls it ``the very singular operator devised by  Hamilton''. Today the nabla operator is almost
universally used in vector calculus, mathematical physics and  and differential geometry, indicating the operation of taking the gradient, divergence, or curl. The exceptions are some continental European textbooks, where these operations are denoted often by abbreviations 
``grad'', ``div'' and ``rot'', instead of nabla.

To conclude this brief review of formulae on coins, I would like to add that formulae are not the only mathematical abstractions which
appear on numismatic materials. I have in my collection a number of coins with designs illustrating geometric concepts,
or presenting  pictorial representations of physical concepts. I plan to describe these in a separate article in the near future.


\providecommand{\href}[2]{#2}\begingroup\raggedright\endgroup

\begin{thebibliography}{1}

\bibitem{Einstein1905}
A.~Einstein, ``Ist die {T}rägheit eines {K}örpers von seinem {E}nergieinhalt
  abhängig?,'' {\em Annalen der Physik} {\bf 18} (1905) 639.

\bibitem{Lev2004}
J.~Leveugle, {\em ``La Relativit\'e, Poincar\'e et Einstein, Planck, Hilbert}.
\newblock Ed. L'Harmattan, 2004.

\bibitem{Euler1740}
L.~Euler, ``De summis serierum reciprocarum,'' {\em Commentarii academiae
  scientiarum Petropolitanae} {\bf 7} (1740) 123--134.

\bibitem{Ostro1845}
M.~V. Ostrogradski {\em Bull. Sci. Acad. Sci. St. Petersburg} {\bf 4} (1845),
  no.~10--11, 145--167.

\bibitem{Leb1901}
P.~Lebedew, ``Untersuchungen \"uber die {D}ruckkr\"afte des {L}ichtes,'' {\em
  Annalen der Physik} {\bf 311} (1901), no.~11, 433--458.

\bibitem{Bogo46}
N.~N. Bogolyubov, ``Kinetic equations,'' {\em JETP} {\bf 16} (1946) 691--702.

\bibitem{mickens1981}
R.~Mickens, {\em An Introduction to Nonlinear Oscillations}.
\newblock Cambridge University Press, 1981.

\bibitem{Banach1922}
S.~Banach, ``Sur les op\'erations dans les ensembles abstraits et leur
  application aux \'equations int\'egrales,'' {\em Fundamenta Mathematicae}
  {\bf 3} (1922) 133--181.

\end{thebibliography}
\end{document}